\documentclass[12pt,leqno]{amsart}  
\usepackage{amsmath,amstext,amsthm,amssymb,amsxtra}
\usepackage[colorlinks,citecolor=red,pagebackref,hypertexnames=false]{hyperref}
\usepackage[author-year]{amsrefs}
\allowdisplaybreaks
\swapnumbers
 
\theoremstyle{plain} 

\newtheorem{theorem}[equation]{Theorem} 

\newtheorem*{sarason}{Characterize Compositions of Paraproducts that are Bounded}
\theoremstyle{definition}

\theoremstyle{remark}

\newtheorem*{tks*}{Acknowledgements}

\numberwithin{equation}{section}

\def\norm#1.#2.{\lVert#1\rVert_{#2}}
\def\Norm#1.#2.{\bigl\lVert#1\bigr\rVert_{#2}}
\def\NOrm#1.#2.{\Bigl\lVert#1\Bigr\rVert_{#2}}
\def\NORm#1.#2.{\biggl\lVert#1\biggr\rVert_{#2}}
\def\NORM#1.#2.{\Biggl\lVert#1\Biggr\rVert_{#2}}

\def\ip#1,#2,{\langle #1,#2\rangle}
\def\Ip#1,#2,{\bigl\langle#1,#2\bigr\rangle}
\def\IP#1,#2,{\Bigl\langle#1,#2\Bigr\rangle}

\def\mid{\,:\,}

\def\abs#1{\lvert#1\rvert}

\def\ABs#1{\Bigl\lvert#1\Bigr\rvert}

\def\XXint#1#2#3{{\setbox0=\hbox{$#1{#2#3}{\int}$}
     \vcenter{\hbox{$#2#3$}}\kern-.5\wd0}}

\def\eqdef{\stackrel{\mathrm{def}}{{}={}}}

\begin{document}
\title[Composition of  Paraproducts] {
Composition of Haar Paraproducts:\\ The Random   Case
}

\author[D. Bilyk ]{Dmitriy Bilyk}

\address{Dmitriy Bilyk  \\
School of Mathematics\\
Georgia Institute of Technology\\
Atlanta GA 30332\\
}

\email {bilyk@math.gatech.edu}

\author[M. T. Lacey]{Michael T. Lacey$^*$}

\address{Michael T. Lacey \\
School of Mathematics\\
Georgia Institute of Technology\\
Atlanta GA 30332\\
}

\thanks{$*.$ Research supported in part by a National 
Science Foundation Grant.  } 
 \email {lacey@math.gatech.edu}

\author[X. C. Li]{Xiaochun Li$^\dagger$}

\address{Xiaochun Li \\
Department of Mathematics\\
University of Illinois, Urbana--Champaign\\
Urbana IL 61801\\
}

\thanks{$^\dagger.$ Research supported in part by a National 
Science Foundation Grant.  } 

\email{xcli@math.uiuc.edu}

\author[B. D.~Wick]{Brett D. Wick$^\ddagger$}

\address{Brett D. Wick \\  Current: 
Fields Institute\\
University of Toronto\\
Toronto, Ontario M5T 3J1 Canada.  \\
Permanent: University of South Carolina\\ Department of Mathematics\\ LeConte College\\ 1523 Greene Street\\ Columbia, SC 29208
}

\thanks{$^\ddagger.$ Research supported in part by a National Science Foundation Grant and the Fields Institute.} 

\email{wick@math.sc.edu}

\maketitle
 
\section{Definitions and Main Theorems} 

We phrase the (difficult) open question which motivates the consideration of this paper. 
Let $ \mathcal D$ be the dyadic grid, and $ \{h_I\mid I\in \mathcal D\}$ 
the $ L ^{2}$ normalized Haar basis, namely 
\begin{equation*}
h_I = \lvert  I\rvert ^{-1/2} \bigl(- \mathbf 1 _{I _{\textup{left}}} 
+ \mathbf 1 _{I _{\textup{right}}}\bigr) 
\end{equation*}
We also use the notation $ h ^{0}_I=h_I$, indicating that the 
Haar function $ h_I$ has integral zero.  Set $ h ^{1}_I \eqdef \abs{ h ^{0}_I}$, 
the subscript $ {} ^{1}$ indicating that $ h^1_I$ has non-zero integral. 

A (classical dyadic) \emph{ paraproduct with symbol $ b$} is one of the operators 
\begin{equation*}
\operatorname B (b, f) = \sum _{I\in \mathcal D} \frac {\ip b, h_I, } 
{\sqrt {\lvert  I\rvert }} \ip f , h ^{\epsilon }_I , h ^{\delta }_I\,. 
\end{equation*}
Here, $ \epsilon ,\delta \in \{0,1\}$, with one of the two being zero and the other 
one.  Then, it is well known that the operator $ \operatorname B (b , \cdot )$ 
is bounded iff the symbol $ b$ is in dyadic $ \operatorname {BMO}$. In particular, 
the following equivalence is a standard part of the literature, and essentially a 
direct consequence of the Carleson Embedding Theorem.  
\begin{equation*}
\norm  \operatorname B (b , \cdot ) . 2 \to 2. 
\simeq 
\sup _{J\in \mathcal D} \Bigl[ \lvert  J\rvert ^{-1} \sum _{I\subset J} \ip b, h_I, ^2
\Bigr]^{1/2} \,. 
\end{equation*}

An outstanding question concerns the \emph{composition} of paraproducts:  Namely, 
if one considers two paraproducts, each of which is potentially unbounded, can one 
meaningfully characterize when the composition is a bounded operator on $ L ^{2}$? 

\begin{sarason}  
Let $  \operatorname B (b , \cdot )$ and $  \operatorname B (b' , \cdot )$ be two possibly
unbounded paraproducts.  
Provide a (non-trivial) characterization of the quantity 
\begin{equation*}
\norm  \operatorname B (b ,  \operatorname B (b' , \cdot ) ). 2\to 2. 
\end{equation*}
(There are distinct versions of this question, depending upon how the 
paraproducts are defined.)
\end{sarason}

We answer a substantially easier variant of this question, in which the 
sum over the dyadic intervals is suitably randomized.  

Motivations for this question arise from formulations of several closely related 
questions.   
(1) Which compositions of Toeplitz operators are bounded? 
\cites{MR511973,MR1395967,MR1934352}.
(2) Which compositions of Hankel operators are bounded? 
\cites{MR898320,MR1395967}
(3) For which pairs of weights $ (u,v)$ is the Hilbert transform 
bounded as a map from $ L ^{2} (du)$ to $ L ^{2}(dv)$?  
The question posed above is a discrete and combinatorial version of 
each of these questions.  For prior work on it, or closely related questions, see 
\cites{MR2097606}.   Each of these questions 
 pose interesting and significant variants, which we leave to the interested 
 reader to pursue in the citations we have provided.  The work of Nazarov, Treil 
 and Volberg on (3), in papers that are both together and separate, are highly recommended. 
 See \cites{MR1685781,MR1945290,MR1428988,MR1428988,MR2019058}.   The paper 
 \cite{MR2174914} studies the role of randomized paraproducts in the theory of product  $
 BMO$. 
 
 There is one variant of these questions that has been completely solved. 
 Eric Sawyer \cite{MR676801} characterized the two weight problem for the Maximal 
 Function.  Also see \cites{MR654182,MR722250}. 

\smallskip

We turn to the main results of this paper. 
Fix a sequence of constants $\mathbf b = \{b_I\mid I\in \mathcal D\} $.  
For choices of $ \epsilon ,\delta \in \{0,1\}$, set 
\begin{equation}\label{e.Pdef}
\operatorname P ^{\epsilon ,\delta }_{\mathbf b}( f) \eqdef \sum _{I\in \mathcal D} 
 {b_I}  \, h ^{\epsilon }_I 
\otimes h ^{\delta }_I 
\end{equation}
There are some comments in order.  First, the typical choice of $ \mathbf b
$ is $ \{\ip b,h_I,/{\sqrt {\abs{ I}}}  \mid I\in \mathcal D\}$, where $ b$ is a $ BMO$ function.  
There is however no function theory in our situation, which is one reason for 
taking $ \mathbf b$ a numerical sequence.  Another reason is that a second 
prominent example, arising from weighted inequalities is $ \mathbf b= \{
\ip b, h ^{1}_I,/{\sqrt {\abs{ I}}} \mid I\in \mathcal D\}$.  We have chosen a normalization 
of $ \mathbf b$ to make it  a `dimensionless' quantity.

It will be important for us that one of our paraproducts be 
permitted to vary over a class of paraproducts.  This class is 
\begin{equation}\label{e.Punc}
\begin{split}
\operatorname P ^{\sigma ,\epsilon ,\delta }_{\mathbf b}( f) \eqdef \sum _{I\in \mathcal D} 
 \sigma _I {b_I}  \, h ^{\epsilon }_I 
\otimes h ^{\delta }_I 
\qquad 
\textup{where}\quad  \sigma =\{\sigma _I\mid I\in \mathcal D\}\in \{-1,+1\} ^{\mathcal D}\,.
\end{split}
\end{equation}
Our primary results concern the 
case where the $\sigma _I $ are independent, identically, uniformly distributed, 
and where the estimates are uniform over the choices of signs.

If a paraproduct does not have a `$ 1$' appearing in its exponents, it is very easy 
to estimate it's norm: 
\begin{equation*}
\norm \operatorname P ^{0,0} _{\mathbf b} .2. 
\simeq 
\sup _I \abs{b_I}\,.
\end{equation*}
The critical case, and the one of interest to us, is when a `$1 $' appears.  
Let us recall the classical result on the boundedness of these paraproducts.  

\begin{theorem}\label{t.classic} We have the equivalence of norms 
\begin{equation}\label{e.classic}
\Norm \operatorname P ^{0,1} _{\mathbf b} .2.
=
\Norm \operatorname P ^{1,0} _{\mathbf b} .2.
\simeq  \sup _{J\in \mathcal D} 
\Bigl[\abs{ J} ^{-1} \sum _{I\subset J} b_I ^2 \abs{ I} 
\Bigr] ^{1/2} 
\end{equation}
\end{theorem}

This is the instance of the Carleson Embedding Theorem.
\begin{equation}\label{e.embed}
\begin{split}
\sum _{I\in \mathcal D} b_I ^2  {\ip f, h^1_I,^2 } 
& \lesssim    \norm \{b_I\}.\textup{Carleson}.^2 \norm f.2.^2 
\\
\norm \{b_I\}.\textup{Carleson}. 
& \eqdef 
\sup _{J\in \mathcal D} 
\Bigl[\abs{ J} ^{-1} \sum _{I\subset J} b_I ^2 \abs{ I} 
\Bigr] ^{1/2} \,. 
\end{split}
\end{equation}
We remark that paraproducts with two $ 1$'s arise in the setting of commutators 
with fractional integrals.  But that is not our theme in this paper.

Our main interest is in the  \emph{composition} of paraproducts, 
with at least one of the two paraproducts having a $ 1$.   The interest here is 
in finding necessary and sufficient conditions for the composition to be bounded, 
permitting the individual paraproducts to be unbounded.  
We concentrate on the random sign case, establishing 
necessary and sufficient conditions which are `natural' extensions of the 
classical result above.

It is useful to set 
notation 
\begin{equation}\label{e.interleave}
\norm  \operatorname  T . {\mathbb E ,p\to p}. 
\eqdef 
\sup _{\norm f.p.=1} \mathbb E \norm \operatorname T f .p. \,. 
\end{equation}
With this definition, it is not necessarily the case $ \norm 
\operatorname  T . {\mathbb E ,p\to p}. \neq
\norm 
\operatorname  T  ^{\ast}. {\mathbb E ,p'\to p'}. $ for conjugate index $p' $.

\begin{theorem}\label{t.main} We have the following equivalences, 
in which $ \sigma $ denotes random choices of signs.
\begin{align}\label{e.10-01}
 \norm \operatorname P ^{1,0} _{\mathbf b}  
\operatorname P ^{\sigma ,0,1} _{\boldsymbol \beta } 
. \mathbb E,2\to2.  
&\simeq 
\norm \{b_I \cdot \beta _I\}. \textup{Carleson}. \,, 
\\ \label{e.01-00}
 \norm \operatorname P ^{0,1} _{\mathbf b}  
\operatorname P ^{\sigma ,0,0} _{\boldsymbol \beta } 
. \mathbb E,2\to2.  
&\simeq \sup _{J} \Bigl[ \frac{\beta _J ^2}{\abs{ J}}  
\sum _{I\;:\;  I\subsetneq J} b _I ^2 \abs{ I} \Bigr] ^{1/2} \,,
\\  \label{e.01-00'}
\norm \operatorname P ^{0,0} _{\mathbf b}  
\operatorname P ^{\sigma ,1,0} _{\boldsymbol \beta } 
. \mathbb E,2\to2.  
&
\simeq \sup _{J} \Bigl[ \frac{\beta  _J ^2} {\abs{ J}}   \sum _{I \;:\;   J\subsetneq I}
b_I ^2 \lvert  I\rvert \Bigr] ^{1/2} \,,
\\ \label{e.01-01} 
\norm \operatorname P ^{0,1} _{\mathbf b} \operatorname P ^{\sigma , 0,1} _{\boldsymbol \beta }. \mathbb E,2\to2.  
&\simeq 
\NOrm \Bigl[ \frac {\beta _J ^2 }  {\abs{ J}} 
\sum _{I\;:\;I\subsetneq J} b_I ^2 \abs{ I} \Bigr] ^{1/2}  . \textup{Carleson}. \,, 
\\
\label{e.01-10}
\begin{split}
 \norm \operatorname P ^{0,1} _{\mathbf b}  
\operatorname P ^{\sigma ,1,0} _{\boldsymbol \beta } 
. \mathbb E,2\to2. 
&\simeq  \bigl\{ \norm \operatorname P ^{\sigma ,0,1} _{\mathbf b}  
\operatorname P ^{0,0} _{\boldsymbol \beta } 
. \mathbb E,2\to2.  +
\norm \operatorname P ^{0,0} _{\mathbf b}  
\operatorname P ^{0,0} _{\boldsymbol \beta } . 2\to2.   
+
\norm \operatorname P ^{0,0} _{\boldsymbol \beta}  
\operatorname P ^{\sigma ,1,0} _{\mathbf b } 
. \mathbb E,2\to2.  \bigr\} \,.
\end{split}
\end{align}
\end{theorem}

In the first equivalence (\ref{e.10-01}) is  
elementary in nature  and we include it only for 
the sake of completeness.  This case dramatically simplifies as it is highly local:
there are no interactions 
between dyadic scales which are widely separated in the the hyperbolic metric.  

The second line, (\ref{e.01-00}) is a characterization in same terms as those of 
Carleson measures. The essential difference is that the non--local:  The coefficients 
$ b_I$ and $ \beta _J$ must be paired even when $ I$ and $ J$ are widely separated in 
scale. 
In this line, it is essential that we have the randomization  
fall  on  the  paraproducts as described above.  Without this term, the right hand side 
is certainly not sufficient for the boundedness of the paraproduct.

Note that for any choices of signs $ \sigma $ and $ \widetilde \sigma $,
 and any $ \phi \in L^2$, we have 
\begin{equation*}
\norm \operatorname P ^{\sigma ,0,1} _{\mathbf b}  
\operatorname P ^{0,0} _{\boldsymbol \beta } \phi 
.2. = 
\norm\operatorname P ^{\widetilde \sigma  ,0,1} _{\mathbf b}  
\operatorname P ^{0,0} _{\boldsymbol \beta } \phi 
.2. 
\end{equation*}
That is, in (\ref{e.01-00}) we need not consider the other randomization 
$ \norm \operatorname P ^{\sigma ,0,1} _{\mathbf b}  
\operatorname P ^{0,0} _{\boldsymbol \beta }
.\mathbb E ,2\to 2.$. 

A similar set of comments apply to (\ref{e.01-00'}).  In particular, the 
characterization here is also of a non-local nature. 

Note that also $ \operatorname P ^{0,0} _{ \mathbf b}
\operatorname P ^{0,0} _{\boldsymbol \beta}$ 
trivially diagonalizes, we we easily have 
\begin{equation*}
\norm 
\operatorname P ^{0,0} _{\mathbf b} \operatorname P ^{0,0} _{\boldsymbol \beta} .2\to 2. 
=\sup _{I} \abs{ b_I \beta _I}\,. 
\end{equation*}

In view of 
(\ref{e.01-00}), the right hand sides of (\ref{e.01-10}) 
can be replaced with an explicit equivalence for the norm in the same spirit as the 
Carleson measure condition.  Again, these conditions are certainly not sufficient 
for the case when one does not average over choices of signs.

There is another instance that arise naturally, through a connection 
with two weight inequalities.  
Let 
\begin{equation*}
\operatorname T _{\sigma } = \sum _{I\in \mathcal D} \sigma _I \, h_I \otimes h_I 
\end{equation*}
by an random Haar multiplier sequence, with $ \sigma _I\in \{\pm 1\}$.  
Let $ \operatorname M _{b} \varphi \eqdef b\cdot \varphi $. 

\begin{theorem}\label{t.Randomntv}  Let $ b, \beta $ be functions, with finitely 
supported Haar expansion. 
We have the equivalence of norms 
\begin{equation}\label{e.Randomntv} 
\begin{split}
\norm \operatorname M_b \operatorname T _{\sigma }\operatorname M_ \beta .\mathbb E ,2\to 2.
&\simeq \sup _{I} \mathbb E \bigl\{  
\norm  \operatorname M_b \operatorname T _{\sigma }\operatorname M_ \beta h ^{1}_I  .2.
+
\norm \operatorname M_\beta  \operatorname T _{\sigma }\operatorname M_ \beta h_I ^{1} .2.
\bigr\} 
\\
& \simeq 
\sup _{I} \frac {\abs{ \ip b,h_I ^1, }} {\abs{ I}} 
\Bigl[\int _I \beta ^2 \; dx  \Bigr] ^{1/2} 
+
\frac {\abs{ \ip \beta ,h_I ^1, }} {\abs{ I}} 
\Bigl[\int _I b ^2 \; dx  \Bigr] ^{1/2} 
\\& \qquad +
\NOrm \frac { \ip b ,h_J,  } {\abs{ J}   } 
\Bigl[ \int_J \beta ^2 \; dx  \Bigr] ^{1/2} . \textup{Carleson}.
\\&\qquad +
\NOrm \frac { \ip \beta ,h_J,  } {\abs{ J}   } 
\Bigl[\int_J b ^2 \; dx \Bigr] ^{1/2}  . \textup{Carleson}.\,.
\end{split}\end{equation}
\end{theorem}

\begin{tks*}
We would like to express our sincere gratitude to the Shanks Foundation at Vanderbilt University.  Their generosity allowed the authors to meet for a weekend at Vanderbilt University to discuss aspects of this question.  Special thanks are due to the organizers of the Thematic Program in Harmonic Analysis at the Fields Institute in Winter-Spring of 2008, as well as the staff of the Institute. It is during that semester that the paper took its final form. 
\end{tks*}

\section{Proof of Theorem~\ref{t.main}} 

The principal way that random signs enters into the the arguments is through 
averaging of all possible signs.  Namely, if $ \sigma _J$ are independent 
identically distributed random variables, taking the values $ \pm1$ with equal 
probability, we have 
\begin{equation*}
\mathbb E \NOrm \sum _{j} \sigma _J g_j .2. ^2 = 
\sum _{j} \norm g_j.2.^2   \,.
\end{equation*}
That is we gain additional orthogonality by averaging over signs.  

To use the random choices of signs on the paraproduct, we need 
random choices of signs indexed by the dyadic intervals
$ \{\sigma  _I \mid I \in \mathcal D\}$, and we write 
\begin{equation*}
\operatorname P ^{\sigma  ,\epsilon ,\delta } _{\mathbf b}
\eqdef 
\sum _{I} \sigma  _I  b_I \, h ^{\epsilon }_I \otimes h ^{\delta }_I
\end{equation*}

\subsection*{Proof of (\ref{e.10-01})}

Observe that the presence of the two inside $ 0$'s diagonalizes the operator. 
\begin{align*}
\operatorname Q & = \operatorname P ^{\sigma,1,0} \operatorname P ^{0,1}
= \sum _{I} \sigma _I \cdot  b_I \beta _I \cdot h ^{1}_I \otimes h ^{1}_I\,.
\end{align*}
That is, the composition is a paraproduct with two $ 1$'s, and symbol is the 
product symbol.  Clearly, random choices of signs  are essential. We have 
\begin{align*}
\mathbb E \norm \operatorname Q \phi .2.^2 
& =  \sum _{I\in \mathcal D} \abs{ b_I \beta _I }^2 \ip \phi ,h_I ^1, ^2 
\end{align*}
And so we can appeal to (\ref{e.embed}) to conclude this case.

\subsection*{Proof of (\ref{e.01-00})}

Observe that 
\begin{equation}\label{e.OOO}
\begin{split}
\operatorname Q & =  \operatorname P ^{0,1} _{\mathbf b}  
\operatorname P ^{\sigma,0,0} _{\boldsymbol \beta }
 = \sum _{\substack{I,J\in \mathcal D\\ I\subsetneq J}}
	\sigma _J \tau _{I,J} \cdot  
	b_I \beta _J  \cdot \tfrac {\sqrt {\abs{ I}}} {\sqrt{\abs{ J}}} \, h_I \otimes h_J 
\end{split}
\end{equation}
where $ \tau _{I,J}$ is the sign of $ h_J$ on the interval $ I$.  Since $ I$ 
is strictly contained in the interval $ J$, this is well defined.

Apply $ \operatorname Q $ to $ h_J$ to see that 
\begin{align*}
\norm \operatorname Q  h_J.2. ^2 
&=\NOrm \sum _{I\subsetneq J} 
\sigma _J \tau _{I,J} \cdot  
	b_I \beta _J  \cdot \tfrac {\sqrt {\abs{ I}}} {\sqrt{\abs{ J}}} \, h_I
	.2. ^2 
\\   
&= \sum _{I\subsetneq J}  \abs{ b_I \beta _J} ^2 \tfrac {\abs{ I}} {\abs{ J}}\,.
 \end{align*}
The operator norm dominates the supremum over $J$ on the right, which proves 
the  lower bound on the operator norm. 

\medskip

Let us assume that 
\begin{equation}\label{e.oo}
\sup _{J}\sum _{I\subsetneq J}  \abs{ b_I \beta _J} ^2 \tfrac {\abs{ I}} {\abs{ J}}\le 1
\end{equation}
and establish an absolute upper bound on the the operator $ \operatorname Q$.  
Apply $ \operatorname Q$ to function $ \phi $ of $ L^2$ norm one. 
\begin{align*}
\mathbb E \norm \operatorname Q \phi .2.^2 
& = \mathbb E  \NOrm  \sum _{\substack{I,J\in \mathcal D\\ I\subsetneq J}}
	\sigma _J \tau _{I,J} \cdot  
	b_I \beta _J  \cdot \tfrac {\sqrt {\abs{ I}}} {\sqrt{\abs{ J}}} \, \ip \phi ,h_J, 
	h_I .2.^2 
\\
&= \sum _{I} \mathbb E \ABs{ \sum _{I\subsetneq J}
\sigma _J \tau _{I,J} \cdot  
	b_I \beta _J  \cdot \tfrac {\sqrt {\abs{ I}}} {\sqrt{\abs{ J}}} \, \ip \phi ,h_J, } ^2 
\\
&=\sum _{\substack{I,J\in \mathcal D\\ I\subsetneq J}} 
(b_I \beta _J) ^2 \tfrac {\abs{ I}} {\abs{ J}} \ip \phi ,h_J,^2 
\\
&\le \norm \phi .2. ^2 \sup _{J} 
\sum _{I\subsetneq J}  \abs{ b_I \beta _J} ^2 \tfrac {\abs{ I}} {\abs{ J}}\,.
\end{align*}
This completes our proof.  

\subsection*{Proof of (\ref{e.01-00'})}

As in the proof of (\ref{e.01-00}), the proof of the lower bound on the 
operator norm  does \emph{not} depend 
upon the averaging over signs. 
Apply $ \operatorname Q= \operatorname P ^{0,0} _{\mathbf b}  
\operatorname P ^{\sigma ,1,0} _{\boldsymbol \beta } $ to $ h_J$ to see that 
\begin{align*}
\norm \operatorname Q h_J.2. ^2 
&=\NOrm \sum _{ I \;:\; J\subsetneq I} 
\sigma _J \tau _{J,I} \cdot  
	b_I \beta _J  \cdot \tfrac {\sqrt {\abs{ I}}} {\sqrt{\abs{ J}}} \, h_I
	.2. ^2 
\\   
&= \sum _{I \;:\; J\subsetneq I}  \abs{ b_I \beta _J} ^2 \frac {\abs{ I}} {\abs{ J}}\,.
 \end{align*}
 Note that the operator 
norm clearly dominates the supremum over $ J$ on the right, and this proves the first  lower bound 
on the operator norm.

Now let us assume that 
\begin{equation}\label{e.OO}
\sup _J \sum _{I \;:\; J\subsetneq I}  \abs{ b_I \beta _J} ^2 \tfrac {\abs{ I}} {\abs{ J}}=1, 
\end{equation}
and let us find an absolute upper bound on operator $ \operatorname Q$ defined 
in (\ref{e.OOO}).  Apply $ \operatorname Q $ to a function $ f$ of $ L^2 $ 
norm one.
\begin{align*}
\mathbb E \norm \operatorname Q f.2.^2 
& = \mathbb E  \NOrm \sum _{\substack{I,J\\ J\subsetneq I }} 
\sigma _J \tau _{I,J} \cdot  
	b_I \beta _J  \cdot \tfrac {\sqrt {\abs{ I}}} {\sqrt{\abs{ J}}} \, \ip f,h_J,\, h_I
	.2.^2 
\\
& =\sum _{I} \mathbb E \ABs{ \sum _{J \;:\; J\subsetneq I}   
\sigma _J \tau _{I,J} \cdot  
	b_I \beta _J  \cdot \tfrac {\sqrt {\abs{ I}}} {\sqrt{\abs{ J}}} \, \ip f,h_J,}^2 
\\
&= \sum _{\substack{I,J\\ J\subsetneq I }} \abs{ b_I \beta _J } ^2 
\tfrac {\abs{ I}} {\abs{ J}} \ip f,h_J,^2 
\\
&\le \norm f.2.^2 \sup _{J} \sum _{I \;:\; J\subsetneq I}  \abs{ b_I \beta _J} ^2 \tfrac {\abs{ I}} {\abs{ J}}\,.
\end{align*}
This completes the proof of the boundedness of this $ \operatorname Q$.

\subsection*{ Proof of (\ref{e.01-01})}

In this proof, we calculate 
\begin{align*}
\mathbb E \norm \operatorname P ^{0,1} _{\mathbf b} \operatorname P ^{\sigma, 0,1} 
_{\boldsymbol \beta } \phi .2.^2 
& = 
\sum _{I} b_I ^2 \abs{ I} \,\mathbb E 
\ABs{ \sum _{J \;:\; I\subsetneq J} \sigma _J \tau _{I,J} \frac {\beta _J} {\sqrt {\abs{ J}}}
 \ip \phi ,h_J ^{1},} ^2 
\\
&=\sum _{J} \ip \phi ,h_J ^{1}, ^2 \cdot \frac {\beta _J ^2 }  {\abs{ J}} 
\sum _{I \;:\; I\subsetneq J} b_I ^2 \abs{ I}\,.
\end{align*}
It is then clear that we have 
\begin{equation*}
\norm \operatorname P ^{0,1} _{\mathbf b} \operatorname P _{\boldsymbol \beta} ^{\sigma, 0,1}
.\mathbb E ,2\to 2. 
\simeq 
\NOrm  \frac {\beta _J  }  {\sqrt{\abs{ J}}} 
\Bigl[ \sum _{I \;:\; I\subsetneq J} b_I ^2 \abs{ I}\Bigr] ^{1/2}  . \textup{Carleson}. 
\end{equation*}

\medskip 

This argument does not shed light on the case when the random choices of signs 
are  imposed on the other paraproduct.  Namely, we obviously have 
\begin{align*}
\norm \operatorname P ^{\sigma,0,1} _{\mathbf b} \operatorname P ^{0,1} 
_{\boldsymbol \beta} \phi .2. ^2 
& = \norm \operatorname P ^{0,1} _{\mathbf b} \operatorname P ^{0,1} 
_{\boldsymbol \beta} \phi .2. ^2 
\\
&=\sum _{I} b_I ^2 \ip h ^{1}_I , \operatorname P ^{0,1} _{\boldsymbol \beta} \phi , ^2\,. 
\end{align*}
That is, the randomization of signs plays no role if placed in this coordinate.

\subsection*{Proof of (\ref{e.01-10})}

Observe that 
\begin{align*}
\mathbb E \norm \operatorname P ^{ 0,1} _{\mathbf b} 
\operatorname P ^{\sigma, 1,0} _{\boldsymbol \beta} \phi .2. ^2 
& = \sum _{I} b_I ^2 \mathbb E 
\ABs{  \sum _{I\cap J\neq \emptyset } \sigma _J
\ip h_I ^{1}, h_J ^{1}, \beta _J \ip \phi , h_J, } ^2
\\
&= \sum _{J} \ip \phi ,h_J, ^2  \beta _J ^2  \sum _{I \;:\;  I\cap J\neq \emptyset } 
b_I ^2 \ip h_I ^{1}, h_J ^{1}, ^2 b _I ^2 \,.
\end{align*}
Therefore, we have 
\begin{align*}
 \norm \operatorname P ^{ 0,1} _{\mathbf b} 
\operatorname P ^{\sigma, 1,0} _{\boldsymbol \beta} .\mathbb E ,2\to 2. ^2 
&\simeq 
\sup _{J}   \beta _J ^2  \sum _{I\mid I\cap J\neq \emptyset } 
b_I ^2 \ip h_I ^{1}, h_J ^{1}, ^2 b _I ^2 
\\
& \simeq \sup _{J} \frac {\beta _J ^2 } {\abs{ J}} 
\sum _{I \;:\; I\subsetneq J} b_I ^{2 } \abs{ I}
+ \beta _J ^2 b_J ^2 + \beta _J ^2 \abs{ J}\sum _{I\mid J\subsetneq I} 
\frac {b_I ^2 } {\abs{ I}} \,. 
\end{align*}
 
By the previous part of the proof, this last supremum, times an absolute constant, 
dominates the sum of the operator norms 
\begin{equation*}
\norm \operatorname P ^{ 0,1} _{\mathbf b} 
\operatorname P ^{\sigma, 0,0} _{\boldsymbol \beta} .\mathbb E,2\to2. ^2 
+ 
\norm \operatorname P ^{ 0,0} _{\mathbf b} 
\operatorname P ^{\sigma, 0,0} _{\boldsymbol \beta} .\mathbb E,2\to2. ^2 
+ 
\norm \operatorname P ^{ 0,1} _{\boldsymbol \beta} 
\operatorname P ^{\sigma, 0,0} _{\mathbf b} .\mathbb E,2\to2. ^2
\end{equation*}
Conversely, by appropriate selection of test function, 
this last term also dominates the supremum, so our proof is complete.

\section{Random Two Weight: Sufficient Direction} 

We consider the upper bound of the operator norms in (\ref{e.Randomntv}) 
by the expressions involving the Carleson measure norm.  To do so, 
we expand the composition $ \operatorname M_b \operatorname T \operatorname M_\beta $ 
as a sum of paraproducts, as considered in the first half of the paper.

Define a paraproduct operator by 
\begin{equation*}
\operatorname P ^{ \beta ,\gamma }_{b,\alpha} 
\eqdef \sum _{I\in \mathcal D} \frac {\ip b,h ^{\alpha }_I,} {{\sqrt {\abs{ I}}}} 
h ^{\beta }_I\otimes h ^ {\gamma }_I 
\end{equation*}
Here, $ \alpha ,\beta ,\gamma \in \{0,1\}$.  Observe that we have the usual expansion 
\begin{align*}
\operatorname M_b \varphi &= 
\Bigl\{ \sum_{I} \ip b,h_I,h_I \Bigr\}
\Bigl\{ \sum_{I} \ip \varphi ,h_I,h_I \Bigr\}
\\
& = \sum _{I} \frac{\ip b,h_I,} {{\sqrt {\abs{ I}}}} 
	\ip \varphi ,h_I,\, h ^{1}_I
\\&\qquad + \sum _{I} \frac{\ip b,h_I,} {{\sqrt {\abs{ I}}}} 
	\ip \varphi ,h_I ^1,\, h ^{1}_I
\\&\qquad + \sum _{I} \frac{\ip b,h_I^1,} {{\sqrt {\abs{ I}}}} 
	\ip \varphi ,h_I ,\, h _I
\\
&= \operatorname P ^{0,1}_{b,0} \varphi +
\operatorname P ^{1,0}_{b,0} \varphi + \operatorname P ^{0,0}_{b,1} \varphi
\end{align*}

There are nine cases in the expansion of the operator
 $ \operatorname M_b \operatorname T \operatorname M_\beta $.  We examine them in turn. 
The symbols associated to the paraproducts reduce to one of four possibilities. 
For the function $ b$ they are 
\begin{equation*}
\mathbf b ^{0} \eqdef \{ \ip b,h ^{0}_I, / \sqrt{\abs{ I}}\}
\,,\, \qquad 
\mathbf b ^{1} \eqdef \{ \ip b,h ^{1}_I, / \sqrt{\abs{ I}}\}\,.
\end{equation*}
 For the function $ \beta $ we have the same two possibilities, 
 for which  we use the notation $ \boldsymbol \beta ^{0}$ and 
$ \boldsymbol \beta ^{1}$.

It is essential to note that if we are forming a composition with 
paraproducts with a $ 0$ in the interior of the composition, then
the random Haar multiplier can be imposed on that zero, 
resulting in a composition of random paraproducts.  Namely, 
\begin{equation*}
\operatorname P ^{ \epsilon _2, 0} _{b,\epsilon _1} 
\operatorname T 
\operatorname P ^{\delta _2,\delta _3} _{\beta,\delta _1} 
=
\operatorname P ^{\sigma,\epsilon _2,0} _{\mathbf b ^{\epsilon _1}} 
\operatorname P ^{\delta _2,\delta _3} _{\boldsymbol \beta ^{\delta _1} }\,.
\end{equation*}
Here, we are using notations from the first half of the paper. 

We discuss the nine separate cases that arise from the expansion of 
\begin{equation} \label{e.9}
\operatorname M_b \operatorname T \operatorname M_\beta = 
\bigl\{  \operatorname P ^{0,0} _{b,1}+ 
 \operatorname P ^{1,0} _{b,0} +\operatorname P ^{0,1} _{b,0}\bigr\}
 \operatorname T 
\bigl\{  \operatorname P ^{0,0} _{\beta,1}+ 
 \operatorname P ^{1,0} _{\beta,0} +\operatorname P ^{0,1} _{\beta,0}\bigr\}
\end{equation}
In each case, we can describe a necessary and sufficient condition 
for the boundedness of that individual term.  These terms together 
provide a sufficient condition for the boundedness of the composition.  
Indeed, the different conditions are listed in  (\ref{e.100-001})
 to (\ref{e.010-001}) and altogether, they are dominated by the 
the right hand side of (\ref{e.Randomntv}).  

\subsubsection*{Case of $00_{1}$---$01_{0}$ and $10_{0}$---$00_{1}$.}

The paraproduct is 
\begin{align*}
\operatorname P ^{0,0}_{b,1} \operatorname T \operatorname P ^{0,1}_{\beta,0} 
&=\operatorname P ^{0,0} _{\mathbf b ^{1}} 
\operatorname P ^{\sigma,0,1} _{\boldsymbol \beta ^{0}}
\end{align*}
This is now a classical paraproduct, and one can estimate it's norm as 
\begin{equation} \label{e.100-001}
\norm 
\operatorname P ^{0,0}_{b,1} \operatorname T \operatorname P ^{0,1}_{\beta,0} 
.2\to2. 
\simeq 
\NOrm \frac {\ip b,h_I ^{1},} {\sqrt {\abs{ I}}} \cdot 
\frac {\ip \beta ,h_I ,} {\sqrt {\abs{ I}}} . \textup{Carleson}. 
\end{equation}

By duality, we see that 
\begin{equation}\label{e.010-100}
\norm 
\operatorname P ^{1,0}_{b,0} \operatorname T \operatorname P ^{0,0}_{\beta,1} 
.2\to2. 
\simeq 
\NOrm \frac {\ip b,h_I ,} {\sqrt {\abs{ I}}} \cdot 
\frac {\ip \beta ,h_I ^{1},} {\sqrt {\abs{ I}}} . \textup{Carleson}. 
\end{equation}

\subsubsection*{Case of $00_{1}$---$00_{1}$.}

The paraproduct is 
\begin{align*}
\operatorname P ^{0,0}_{b,1} \operatorname T \operatorname P ^{0,0}_{\beta,1} 
&=\operatorname P ^{\sigma,0,0} _{\mathbf b ^{1}} 
\operatorname P ^{0,0} _{\boldsymbol \beta ^{1}}
\end{align*}
This is a composition with all zeros, which immediately diagonalizes to give 
\begin{equation}\label{e.100-100}
\norm \operatorname P ^{0,0}_{b,1} \operatorname T \operatorname P ^{0,0}_{\beta,1} 
.2\to2. 
=  \sup _{I}  \frac { \abs{ \ip b, h_I ^{1}, \ip \beta ,h_I ^{1},}} {\abs{I} } 
\end{equation}

\subsubsection*{Case of  $00_{1}$---$10_{0}$ and $01_{0}$---$00_{1}$. }

The paraproduct in the case of $00_{1}$---$10_{0}$  is 
\begin{align*}
\operatorname P ^{0,0}_{b,1} \operatorname T \operatorname P ^{1,0}_{\beta,0} 
&=\operatorname P ^{\sigma,0,0} _{\mathbf b ^{1}} 
\operatorname P ^{1,0} _{\boldsymbol \beta ^{0}}
\end{align*}
This is a paraproduct of the type considered in (\ref{e.01-00}) (with the 
roles of $ b$ and $ \beta $ reversed.) Thus, we have 
\begin{equation}\label{e.100-010}
\begin{split}
\norm \operatorname P ^{0,0}_{b,1} \operatorname T \operatorname P ^{1,0}_{\beta,0}  
.2\to2. ^2 
& \simeq 
	\sup _{J} \frac{{ \ip b,h_J ^1, ^2 }} {\abs{ J} ^2 } 
	\sum _{I\subsetneq J} \ip \beta ,h_I, ^2 
\\
& \simeq 
	\sup _{I} { \ip \beta ,h_I, ^2 } 
	\sum _{I\subsetneq J}  
	 \frac{{ \ip b,h_J ^1, ^2 }} {\abs{ J} ^2 }
\end{split}
\end{equation}

By duality, in the case of $01_{0}$---$00_{1}$ we have 
\begin{equation}\label{e.001-100}
\begin{split}
\norm \operatorname P ^{0,1}_{b,0} \operatorname T \operatorname P ^{0,0}_{\beta,1}  
.2\to2. ^2 
& \simeq 
	\sup _{J} \frac{  \ip \beta,h_J ^1, ^2 } {\abs{ J} ^2 } 
	\sum _{I\subset J} \ip b ,h_I, ^2 
\\
& \simeq 
	\sup _{I} { \ip b ,h_I, ^2 } 
	\sum _{I\subsetneq J}  
	 \frac{{ \ip \beta,h_J ^1, ^2 }} {\abs{ J} ^2 }
\end{split}
\end{equation}

\subsubsection*{Case of $01_{0}$---$10_{0}$}

This case does not appeal to the first half of our paper, as the 
random Haar multiplier cannot be absorbed into either paraproduct.  We estimate 
\begin{align*}
\mathbb E \norm \operatorname P ^{0,1}_{b,0} \operatorname T \operatorname P ^{1,0}_{\beta,0} \phi .2.^2 
& 
=\sum _{I}  \frac {\ip b,h_I , ^2 } {\abs{ I} } 
 \mathbb E \Ip h_I ^{1} , \operatorname T 
 	\operatorname  P ^{1,0}_{\beta,0}  \phi , ^2 
\\
&= \sum _{I} {\ip b,h_I, ^2 }   
\sum _{J\mid I\subsetneq J} \frac {\ip \beta ,h _{J}, ^2 } {\abs{ J }} 
	\ip \phi ,h_J, ^2 
\end{align*}
It is then clear that we have 
\begin{equation}\label{e.001-010} 
\norm \operatorname P ^{0,1}_{b,0} \operatorname T \operatorname P ^{1,0}_{\beta,0} .2\to 2. ^2  
\simeq 
\sup _{J} \frac {\ip \beta ,h _{J}, ^2 } {\abs{ J }} 
\sum _{I \;:\; I\subsetneq J} {\ip b,h_I, ^2 } 
\end{equation}

\subsubsection*{The Case of $01_{0}$---$01_{0}$ and  $10_{0}$---$10_{0}$. }

The case of $01_{0}$---$01_{0}$ is 
\begin{equation*}
\operatorname P ^{0,1} _{b,0} \operatorname T \operatorname P ^{0,1} _{\beta ,0} 
= \operatorname P ^{0,1} _{\mathbf b ^{0}} 
	\operatorname P ^{\sigma,0,1} _{\boldsymbol \beta ^{0}}.
\end{equation*}
Thus, we can appeal to   (\ref{e.01-01}) to conclude that 
\begin{equation}\label{e.010-010}
\norm 
\operatorname P ^{0,1} _{b,0} \operatorname T \operatorname P ^{0,1} _{\beta ,0} .2.   
\simeq 
\NOrm \frac { \ip \beta ,h_J,  } {\abs{ J}   } 
\Bigl[\sum _{I \;:\; I\subsetneq J} \ip b ,h_I, ^2 \Bigr] ^{1/2}  . \textup{Carleson}. \,.
\end{equation}

The case of  $10_{0}$---$10_{0}$ is dual, and yields 
\begin{equation}\label{e.001-001}
\norm 
\operatorname P ^{1,0} _{b,0} \operatorname T \operatorname P ^{1,0} _{\beta ,0} .2. ^2 
\simeq 
\NOrm \frac { \ip b ,h_J,  } {\abs{ J}   } 
\Bigl[\sum _{I\subsetneq J} \ip \beta  ,h_I, ^2  \Bigr] ^{1/2} . \textup{Carleson}. \,.
\end{equation}

\subsubsection*{The Case of $10_{0}$---$01_{0}$. }

This case reduces to that of (\ref{e.10-01}). 
\begin{equation}\label{e.010-001}
\begin{split}
\operatorname P ^{1,0} _{b,0} \operatorname T \operatorname P ^{0,1} _{\beta ,0} 
&= \operatorname P ^{\sigma,1,0} _{\mathbf b ^{0}} 
	\operatorname P ^{0,1} _{\boldsymbol \beta ^{0}}.
\\
& \simeq 
\NOrm \frac {\ip b ,h_I, \ip \beta ,h_I, } {\abs{ I}} . \textup{Carleson}. \,.
\end{split}
\end{equation}

\section{Random Two Weight: The Necessary Direction} 

We show that the operator norms on the left in (\ref{e.Randomntv}) provide an upper 
bound on the expressions involving Carleson measures.  

Let us begin by applying the operators $ \operatorname M_b \operatorname  T 
\operatorname M_ \beta $ to the functions $ h ^{1}_J$. We have 
\begin{align*}
\mathbb E \norm \operatorname M_b \operatorname  T 
\operatorname M_ \beta .2\to 2. ^2 
&\ge 
\mathbb E \norm \operatorname M_b \operatorname  T 
\operatorname M_ \beta  h ^{1}_J .2. ^2 
\\
& \ge \frac 1 {\abs{ J}} \sum _{I\subset J} \norm b  \cdot h_I \cdot  \ip \beta ,h_I, .2.^2 
\\
& = \frac 1 {\abs{ J}} \sum _{I\subset J} \frac{\ip \beta ,h_I, ^2 } {\abs{ I}} 
\int _I b ^2 \; dx 
\\
&=
\frac 1 {\abs{ J}} \sum _{I\subset J} \frac{\ip \beta ,h_I, ^2 } {\abs{ I}}
\Bigl\{  \ip b ,h_J ^1, ^2 +\sum _{K\subsetneq J } \ip b,h_J,^2 \Bigr\}
\end{align*}
This inequality, and the dual assertion, 
proves that the operator norm bounds the two terms involving the 
Carleson measure on the right in (\ref{e.Randomntv}). 

It remains to observe that we have 
\begin{align*}
\mathbb E \norm \operatorname M_b \operatorname  T 
\operatorname M_ \beta .2\to 2. 
&\ge 
\mathbb E \norm \operatorname M_b \operatorname  T 
\operatorname M_ \beta  h ^{0}_J .2. ^2 
\\ 
&\ge \frac { \int _J b ^2 \; dx \abs{ \ip \beta ,h_J ^{0} ,}} {\abs{ J} ^2 }
\end{align*}
This, and the dual inequality, completes the proof of the operator bound bounding the the 
right hand side of (\ref{e.Randomntv}).

\end{document}